\documentclass{article}
\usepackage[usenames,dvipsnames]{xcolor}
\usepackage{amsmath, amssymb, amsfonts, color, stmaryrd}

\newcommand{\be}{\begin{equation}}
\newcommand{\ee}{\end{equation}}

\newcommand{\ba}{\begin{align*}}
\newcommand{\eal}{\end{align*}}

\newcommand{\RRM}{\mathbb{R}}

\newcommand{\cll}{{\cal L}}

\newcommand{\clg}{{\cal G}}
\newcommand{\cle}{{\mathcal E}}

\newcommand{\clz}{{\cal Z}}

\newcommand{\clr}{{\cal R}}

\newcommand{\vdiv}{{\mathop{\mbox{\rm div\,}}}}

\newcommand{\clp}{{\cal P}}
\newcommand{\cla}{{\cal A}}

\newcommand{\clb}{{\cal B}}

\newcommand{\supp}{\mathop{\mbox{\rm supp}}}

\newcommand{\bea}{\begin{eqnarray*}}
\newcommand{\eea}{\end{eqnarray*}}
\newcommand{\eps}{\varepsilon}

\newcommand{\Ad}{\mbox{\rm Ad}}

\newcommand{\bba}{{\mathbb A}}
\newcommand{\bbb}{{\mathbb B}}
\newcommand{\bbe}{{\mathbb E}}
\newcommand{\bbf}{{\mathbb F}}
\newcommand{\bbz}{{\mathbb Z}}
\newcommand{\bbl}{{\mathbb L}}
\newcommand{\clc}{{\cal C}}
\newcommand{\clx}{{\cal X}}
\newcommand{\cly}{{\cal Y}}

\newtheorem{theorem}{Theorem}[section]
\newtheorem{lemma}[theorem]{Lemma}

\newtheorem{definition}[theorem]{Definition}

\newenvironment{keywords}%
   {\begin{trivlist}\item[]{\bfseries\sffamily Keywords:}\ }
   {\end{trivlist}}

\newenvironment{MSC classification}%
   {\begin{trivlist}\item[]{\bfseries\sffamily MSC classification:}\ }
   {\end{trivlist}}

\numberwithin{equation}{section}

\begin{document}

\title{A moving boundary problem for the Stokes equations involving  osmosis:
variational modelling and short-time well-posedness}
\author{Friedrich Lippoth \\ \small{Institute of Applied Mathematics, Leibniz 
University Hannover, Welfengarten 1,} \\ \small{ D-30167 Hannover, Germany }
\\ \small{e-mail: \texttt{lippoth@ifam.uni-hannover.de}} \\
\\Mark A. Peletier\\ 
Georg Prokert\\
\small{Faculty of Mathematics and Computer
Science, TU Eindhoven} \\ \small{P.O. Box 513 5600 MB Eindhoven, the
Netherlands} \\
\small{e-mail: \texttt{m.a.peletier@tue.nl, g.prokert@tue.nl}}
}

\date{}
\maketitle

\begin{center}
\end{center}

\begin{abstract}
Within the framework of variational modelling we derive  a one-phase moving
boundary problem describing the motion of a semipermeable membrane enclosing a
viscous liquid, driven by osmotic pressure and surface tension of the membrane.
For this problem we prove the existence of classical solutions for a short time.
\end{abstract}

\begin{keywords}
variational modelling, Stokes equations, osmosis,  moving boundary problem,
maximal continuous regularity
\end{keywords}

\begin{MSC classification}
35R37, 35K55, 76M30, 80M30 
\end{MSC classification}

\section{Introduction} \label{s1}
It is the aim of the present paper to introduce and discuss a model for the
motion of a closed membrane in a liquid, taking into account simultaneously the
following effects:

\begin{enumerate}
\item surface tension forces of the membrane,
\item diffusion of a solute in the liquid,
\item (quasistationary) viscous motion of the liquid,
\item osmotic pressure difference across the membrane,
\item resistance to liquid motion through the membrane.
\end{enumerate}

To avoid additional difficulties, our attention in this paper is restricted to a
one-phase problem, i.e. we assume that the liquid outside the membrane has
negligible viscosity and does not contain any solute.

We remark that quite a number of well studied moving boundary problems are
contained as special or limit cases in the above general setting. For example,
if 1. and 3. are the only forces taken into account and the membrane is
 impermeable, so-called quasistationary Stokes flow (driven by surface
tension) arises. If
only 1., 2., and 4. are considered, an osmosis model without liquid
motion occurs. In the limit of fast diffusion, this problem
(formally) yields a mean curvature flow problem with a nonlocal
lower order term. For details and references on this, we refer to Section
\ref{modlim}.

It is reasonable to demand that a model for the full problem should recover
these special and limiting cases. The crucial observation here is that both the
osmosis problem and the Stokes flow problem
 have a variational structure, i.e. they can be
interpreted as gradient flows with respect to certain energies and dissipation
functionals. This yields a straightforward approach to the full problem:
We use linear combinations for both energy and dissipation  and
assume that the evolution we are interested in is a generalized gradient flow
with respect
to these. More precisely, an additional dissipation term is introduced to model
(finite) resistance of the membrane against solute flux. 

We will use this approach to formally derive a moving boundary PDE system 
that describes our full problem. At the moment it seems quite challenging  to
use its variational
structure for deriving existence results under weak smoothness assumptions. 
Therefore we do not pursue this approach here. Instead, we show
short-time well-posedness of the PDE system in a classical setting by
transformation to a fixed domain and applying maximal-regularity results and a
contraction argument to the resulting nonlinear parabolic Cauchy problem in a
product of suitable function spaces, in a fashion oriented at \cite{lipr1}. 

This paper is organized as follows: In Section \ref{modlim} we introduce the 
model and discuss several limit cases. In Section \ref{pre} we collect some
necessary technical material such as a regularity result in the scale of little
H\"older spaces for the Stokes equations. We define the notion of a classical
solution of the full problem and state our main result (Theorem \ref{mt}).
Section \ref{trafo} contains its proof.

\section{Modelling and limit cases}\label{modlim}

%
%
%
%
%
%
%
%
%
%
%
%
%
%
%
%
The following considerations and calculations that motivate our model are 
formal,  i.e. we will
assume sufficient smoothness and do not strictly define the phase manifold, its
tangent space etc. We follow the methodology of \emph{variational modelling};  
see~\cite{PeletierVarMod} for a more detailed introduction and overview. 
\subsection{Variational modelling: basic approach and preliminary examples} \label{ss21}
\paragraph{The general setting.} Defining a variational model comprises the 
following choices:
\begin{itemize}
\item  a set $\clz$ of states with the structure of a 
differentiable (possibly infinite-dimensional) manifold (the {\em state 
manifold});
\item  a differentiable {\em energy functional} 
$\cle:\;\clz\longrightarrow\RRM$; 
\item  a {\em process space} $\clp_z$ at all $z\in\clz$; (In fact, 
it makes sense to consider $\clp:=\bigcup_{z\in\clz} \clp_z$ as a vector bundle 
over $\clz$.)
\item  a {\em dissipation functional} 
$\Psi_z:\;\clp_z\longrightarrow\RRM$,
\item  a linear (bundle) map 
$\Pi_z:\;\clp_z\longrightarrow T_z\clz$, where $T_z\clz$ denotes the tangent 
space to $\clz$ at $z$. We will refer to $\Pi_z$ as {\em process map.}
\end{itemize}
With these choices, the model is given as a dynamical system on $\clz$ is
defined by
\begin{equation}\label{varmin1}
\dot z=\Pi_z w^\ast,
\end{equation}
where $w^\ast$ is the solution to the minimization problem
\begin{equation}\label{varmin2}
\Psi_z(w)+\cle'(z)[\Pi_z w]\;\longrightarrow\;\min,\qquad w\in\clp_z.
\end{equation}
Here $\cle'(z)[s]$ is the Fr\'echet derivative of $\cle$ at $z$ applied to the 
tangent vector $s$. The minimizer $w^\ast$ can be considered as the 
``actual process chosen by the system in state $z$'' in the trade-off between
diminishing the 
energy and minimizing dissipation at $z\in\clz$.
We assume here that the minimization problem (\ref{varmin2}) is uniquely 
solvable, a property that hinges on strict convexity and coercivity in 
appropriate norms of the potential $\Psi_z$. In various other applications, such 
as the rate-independent systems that appear in fracture, plasticity, and 
hysteresis~\cite{Mielke05a}, these properties fail;  we refrain 
from giving details on this as our interest is restricted to a situation in 
which the above framework is sufficient. 

It is an advantage of this modelling approach that in many applications it 
results in modelling choices that are independent from each other, relatively 
straightforward, and can be justified from physical principles. By contrast,
the resulting systems of differential equations might be less 
transparent. 

Observe, furthermore, that the approach described here contains the concept of 
a {\em gradient flow} as a special case: If $(\clz,g)$ is a Riemannian 
manifold, $\clp_z=T_z \clz$, $\Psi_z(w)=\frac{1}{2}g_z(w,w)$, and $\Pi_z$ is the 
identity, then (\ref{varmin1}) defines the gradient flow on $(\clz,g)$ induced 
by $\cle$. On the other hand, if the dissipation functional $\Psi_z$ is 
quadratic (i.e. $\Psi_z(w) = \Phi_z(w,w)$) and positive definite, and if the process 
map is surjective, any evolution of type (\ref{varmin1}, \ref{varmin2}) can be 
considered as a gradient flow with respect to $\cle$ and the Riemannian metric $g_z$  
constructed as follows: let $Q_z$ be the $\Psi_z$-orthogonal projection along 
$\mbox{ker} \Pi_z$. Observe that $Q_z$ is constant along fibres 
$\Pi_z^{-1}(v)$, $v \in T_z \clz$. Using this, it is straight forward to check that the bilinear map 
\[
g_z(v_1,v_2) := \Phi_z(Q_z w_1, Q_z w_2), \qquad \Pi_z w_i = v_i \;\; (i = 1,2), 
\]
is a well-defined Riemannian metric on $T_z \clz$ that has the desired property.

Before we are going to describe the model we are interested in we will 
informally illustrate the approach of variational modelling in two related but 
simpler, paradigmatic and well understood problems. This will also provide a 
motivation for the choices we are going to make for the osmosis model. 

\paragraph{Diffusion~\cite[Sec.~5.5]{PeletierVarMod}.} 
 To model linear diffusion on $\RRM^N$, 
the following choices are possible within the above framework. $\clz$ is a 
suitable set of scalar nonnegative functions on $\RRM^N$ representing 
concentrations. The functional $\cle$ is given by the entropy
\[\cle(c):=\gamma\int_{\RRM^N}c\ln c\,dx,\qquad\mbox{$\gamma>0$ fixed, 
$c\in\clz$.}\]
(We remark that $c$ is
understood to be dimensionless, i.e. a ratio with respect to a fixed reference
concentration $c^0$. For simple, approximately spherical particles, we
have $\gamma =
RTc^0$, where $R$ is the universal gas constant, $T$ is absolute temperature,
and $c^0$ is the same normalization concentration, see
e.g.~\cite[Ch.~4]{PeletierVarMod}).
The process space consists of vector fields on $\RRM^N$ representing the mass 
flux, with dissipation functional
\[\Psi_c(f):=\frac{\eta_1}{2}\int_{\RRM^N}\frac{|f|^2}{c}\,dx,\qquad f\in
\clp_c.\]\
(This is $+\infty$ whenever $c=0$ and $f\neq 0$ on a set of nonzero measure, so
for the actual process, $f=0$ is enforced where $c=0$.)
The constant $\eta_1$ in this expression is an inverse mobility of the solute.
 For dilute,
approximately spherical solute  particles
it should be taken equal to $6\pi r \nu$,
where $r$ is the particle radius, and $\nu$ the dynamic viscosity of the
solvent~\cite[Ch.~5]{PeletierVarMod}.

Naturally, the tangent space $T_c(\clz)$ has to be interpreted as the set of 
the ``local concentration changes''. Accordingly, to encode mass conservation 
we choose 
\[\Pi_cf:=-\vdiv f,\qquad f \in \clp_c.\]
Solving (\ref{varmin2}) with these choices yields $\eta_1f^\ast=\gamma\nabla 
c$ and from (\ref{varmin1}) we get 
\[c_t=\frac{\gamma}{\eta_1} \Delta c\quad\mbox{on $\RRM^N$}.\]

\paragraph{Free-boundary Stokes flow driven by curvature.}  To 
describe the motion of a liquid drop that is deformed by  surface tension 
forces via ``creeping flow'', it is natural to choose $\clz$ to be the 
(infinite\-dimensional) manifold of (simply connected) domains $\Omega$, 
representing the 
drop shapes. The energy is the corresponding surface measure:
\[\cle(\Omega):=\alpha|\partial\Omega|=\alpha\int_{\partial\Omega}\,d\sigma,\]
where $\alpha>0$ is the surface energy density.
The tangent space $T_\Omega\clz$ can be represented by functions 
$V_n:\;\partial\Omega\longrightarrow\RRM$ that play the role of the normal 
velocities of the boundary of a moving domain $t\mapsto\Omega(t)$. 

Again it is straightforward to choose a suitable process space, namely, the 
space of all divergence\-free vector fields in $\Omega$, with dissipation 
caused by inner friction of the liquid:
\[\Psi_\Omega(v):=\frac{\eta_2}{2}\int_\Omega|\eps(v)|^2\,dx,\qquad
\eps(v)=\tfrac{1}{2}(\nabla v+\left(\nabla v)^T\right),\qquad v\in\clp_\Omega.\]
These choices are implied by the assumptions that the liquid is Newtonian and 
incompressible, and that its mass is preserved. The constant $\eta_2$ is the
shear viscosity. The mapping $\Pi_\Omega$ 
is chosen to represent the kinematic boundary condition, i.e. the
assumption that the 
boundary of the drop moves along with the liquid particles that constitute it:
\[\Pi_\Omega(v):=v|_{\partial\Omega}\cdot n,\]
where $n$ is the exterior unit normal to $\Omega$. 

A straightforward calculation shows that (\ref{varmin1},\ref{varmin2}) now 
yield 
 the moving boundary problem 
\[V_n=u\cdot n\quad\mbox{on $\partial\Omega$},\]
where $u$ solves the Stokes system
\bea
\begin{array}{rcll}
\eta\Delta u-\nabla p&=&0&\mbox{ in $\Omega$,}\\
\vdiv u&=&0&\mbox{ in $\Omega$,}\\
\tau(u,p)n&=&\alpha Hn&\mbox{ on $\partial\Omega$,}
\end{array}
\eea
where $p$ is the hydrodynamic pressure, occurring here as a Lagrange multiplier 
coresponding to the incompressibility condition, $H$ is the sum of the 
principal curvatures of $\partial\Omega$, and 
\[\tau(u,p)=\eta_2\eps(u)-pI\]
is the hydrodynamic stress tensor. We have to remark here that the Stokes 
system defines $u$ only up to rigid body velocities  since these are in the 
kernel 
of $\eps$, see below. This moving boundary problem (as well as closely related
ones) has been discussed extensively in the literature, see e.g.~\cite{ant,
espr, gupr, hopper, sol}.

In both problems discussed here, the following observations can be made: 
\begin{itemize}
\item Quadratic dissipation functionals correspond to linear constitutive 
relations (Fick's law and Newtonian stress--strain relation, respectively).
\item In these cases and many others, the maps $\Pi_z$ encode balance 
laws (mass balance in the first 
example, ``boundary conservation'' in the second). 
\end{itemize}

\subsection{Variational modelling of the Stokes-Osmosis problem}
In the process of variational modelling, the above examples will play a guiding
role. We will need to deal with some  additional aspects: mass conservation of
the solvent in the presence of a moving boundary, dissipation by the osmotic
process, and diffusion in a moving solvent.

We start by listing our modelling assumptions.
\begin{itemize}
\item The solvent is incompressible and moves with velocity field $u$.
Mass conservation of the solvent then implies 
\begin{equation}
\label{cond:incomp}
\vdiv u=0\quad\mbox{in $\Omega$.}
\end{equation}
\item The solute moves according to a flux field $f$ so that mass
conservation of the solute is expressed by
\begin{equation}
\label{rel:c_t-f}
c_t+\vdiv f=0\quad\mbox{in $\Omega$.}
\end{equation}
\item The membrane is impermeable to the solute. Together with mass
conservation this implies the boundary condition
\begin{equation}
\label{rel:f-Vn}
f\cdot n=cV_n\quad\mbox{on $\partial\Omega$.}
\end{equation}
\item The solvent inside the cell is viscous with a linear dependence of
the strain rate on the stress. 
\item The solute motion is governed by convection along $u$ and diffusion
through the solvent. The diffusion obeys Fick's law. 
\item There is finite resistance by the membrane to solvent moving through it,
proportional to the solvent flux.
\end{itemize}

These assumptions obviously oversimplify the physics of any ``real'' membrane
enclosing a liquid in motion. No further mechanical properties except
resistance to area growth from normal displacement are taken into account; in
particular, there is no resistance to tangential stretching (or one
has to assume completely frictionless slipping of the liquid in tangential
direction along the membrane). Including these effects, in the vein of e.g.
\cite{ccs} seems to be interesting but highly nontrivial and has to remain
outside the scope of this paper.

Our model is encoded in the following choices within the variational
framework described above. 
\paragraph{State manifold.} Since we consider a coupled problem involving
diffusion inside a moving domain, we choose the manifold of  pairs $(\Omega,c)$
where $\Omega$ is a bounded domain in $\RRM^N$ and $c$ is a nonnegative solute
concentration such that $\supp c\subset\bar\Omega$. Accordingly its tangent
space consists of pairs $(V_n, c_t)$ where $V_n$ is  has the same meaning as
above and $c_t$ is a scalar function in $\Omega$ representing concentration
changes. 
\paragraph{Energy functional.} Keeping in mind the examples above, we define
\[\cle(\Omega,c)=\gamma\int_\Omega c\ln c\,dx+\alpha|\partial\Omega|.\]
with appropriate positive constants $\alpha$ and $\gamma$. This includes
diffusion and surface tension as driving mechanisms of the evolution.
\paragraph{Process space.} The processes that cause energy dissipation here are
solvent motion, solute flux, and passage of the solvent through the membrane.
As in the examples above, the former two are described by a velocity field $u$
and a flux field $f$. Since the solvent flux through the membrane is given by 
$u|_{\partial\Omega}\cdot n-V_n$  and in view of (\ref{rel:f-Vn}) it makes
sense to consider $V_n$ as a process component as well. We therefore  
choose (cf. (\ref{cond:incomp},\ref{rel:f-Vn}))
\[\clp_{(\Omega,c)}=\{(u,f,V_n)\,|\,\vdiv u=0\mbox{ in $\Omega$, }
f\cdot n=cV_n\mbox{ on $\partial\Omega$}\}.\]

\paragraph{Process map.} 
In view of (\ref{rel:c_t-f}) it is now straightforward to define
\[\Pi_{(\Omega,c)}(u,f,V_n)=(V_n,-\vdiv f).\]

\paragraph{Dissipation.} Concerning  dissipation by  the motion of the solutes, 
we
have to consider now the flux relative to the flux $cu$ arising from pure
transport by the solvent. The dissipation by inner friction is the same as
described in the example above. Additionally, we  model the resistance
that solvent particles have to overcome when they cross the membrane.
Since we assume a linear constitutive relation here as well,
we find 
\[\Psi_{(\Omega,c)}(u,f,V_n):=
\frac{\eta_1}{2}\int_\Omega\frac{| f-cu|^2}{c}\,dx +
\frac{\eta_2}{2}\int_\Omega|\eps(u)|^2\,dx+
\frac{\eta_3}{2}\int_{\partial\Omega}(u\cdot n-V_n)^2\,d\sigma.\]
The constants $\eta_1$ and $\eta_2$ have the same meaning as in the
introductory examples above, and $\eta_3$ can be interpreted as the inverse
permeability of the membrane. 

These choices complete the modelling of our problem, in the  sense that the 
model is fully specified by them. We next derive the evolution equations for the state~$(\Omega,c)$.

\paragraph{Remark.} \begin{itemize} \item[i)] Observe that for $c>0$, the normal velocity $V_n$ is
determined by (\ref{rel:f-Vn}). This suggests the alternative (but, in this
case, equivalent) choices
\begin{align*}
\tilde\clp_{(\Omega,c)}&=\{(u,f)\,|\,\vdiv u=0\mbox{ in $\Omega$}\},\\
\tilde\Pi_{(\Omega,c)}&=\left({\textstyle\frac{1}{c}}f|_{\partial\Omega}\cdot
n,-\vdiv f\right).
\end{align*}
This appears to be more elegant and straightforward because the process space
is smaller and the process map is completely dictated by solute mass
conservation. However, we refrain from this as the restriction to everywhere
strictly positive concentrations is both unnatural and unnecessary in our
analysis of the resulting moving boundary problem. 

\item[ii)] If we restrict the manifold of states to any 
submanifold  consisting of pairs $(\Omega,c)$ that satisfy
\[\int_\Omega c\,dx=M_0\]
with $M_0$ being a fixed total solute mass
then the process maps 
$\Pi_{(\Omega,c)}$  are easily seen to be surjective. 
\end{itemize}

We turn to the minimization problem (\ref{varmin2}) and observe first that
\[\cle'(\Omega,c)[\Pi_{(\Omega,c)}(u,f,V_n)]
=\gamma\int_\Omega\frac{\nabla
c}{c}\cdot f\,dx-\int_{\partial\Omega}(\alpha H+\gamma c)V_n\,d\sigma.\]
As usual, we   account for the incompressibility condition
(\ref{cond:incomp}) by introducing a Lagrange multiplier $q$, which physically
represents the hydrodynamic pressure. Thus the stationarity conditions are the
vanishing of the first variation of
\[L(u,f,V_n,q):=\Psi(u,f,V_n)+
\gamma\int_\Omega\frac{\nabla
c}{c}\cdot f\,dx-\int_{\partial\Omega}(\alpha
H+\gamma c)V_n\,d\sigma-\int_\Omega
q\,\vdiv u\,dx\]
with respect to all variations $(\tilde u,\tilde f,\tilde V_n)$ that satisfy
$\tilde f\cdot n=c\tilde V_n$  on $\partial\Omega$. 

 Explicitly, this means
\begin{align}
\int_{\partial\Omega}(-\eta_3(u\cdot n-V_n)-\alpha H-\gamma
c)\tilde V_n\,d\sigma&=0,\nonumber\\
-\eta_1\int_\Omega(f-cu)\cdot\tilde
u\,dx+\eta_2\int_\Omega\eps(u):\eps(\tilde u)\,dx-\eta_3\int_{\partial\Omega}
(V_n-u\cdot n)\tilde u\cdot n\,d\sigma-\int_\Omega q\vdiv\tilde
u\,dx&\nonumber\\
=\int_\Omega(-\eta_2\Delta u +\nabla q-\eta_1(f-cu)) \cdot\tilde u\,dx 
+\int_{\partial\Omega}(\tau(u,q)n-\eta_3(V_n-u\cdot n)n)\cdot\tilde
u\,d\sigma&=0,\nonumber\\
\int_\Omega\frac{\eta_1(f-cu)+\gamma\nabla c}{c}\cdot\tilde
f\,dx&=0,\label{fvar}
\end{align}
 where
\[ \tau(u,q)=\eta_2\eps(u)-qI\]
is the (hydrodynamic) stress tensor in the solvent as before. 

Gathering all equations and eliminating $f$ by means of (\ref{fvar})$_3$ we
obtain the moving boundary problem
\be\label{mbp0}
\left.\begin{array}{rclll}
-\eta_2\Delta u+\nabla(q+\gamma c)&=&0&\mbox{ in $\Omega(t)$,} & t > 0,\\
\vdiv u&=&0&\mbox{ in $\Omega(t)$,} & t > 0,\\
\partial_t c-\frac{\gamma}{\eta_1}\Delta c+\nabla c\cdot u&=&0&\mbox{ in 
$\Omega(t)$,} & t > 0,\\
\tau(u,q+\gamma c)n-\alpha Hn&=&0&\mbox{ on $\partial \Omega(t)$,} & t > 0,\\
\alpha H+\gamma c+\eta_3(u\cdot n-V_n)&=&0&\mbox{ on $\partial \Omega(t)$,} & t 
> 0,\\
-\frac{\gamma}{\eta_1}\partial_nc+cu\cdot n-cV_n&=&0&\mbox{ on $\partial 
\Omega(t)$,} & t > 0.
\end{array}\right\}
\ee
To enforce uniqueness of the solution, in our one-phase setting one has to
exclude rigid body motions. For this, we shall additionally demand that $u$ is
$L^2(\Omega)$-orthogonal to the space of rigid body velocities on $\RRM^N$.
Of course, the system has to be complemented by initial conditions for $c$ and 
$\Omega$. 

The model proposed here contains a number of more\-or\-less well-studied moving
boundary problems as (formal) limit cases: 

\begin{itemize}
\item {\bf Osmosis in a resting solvent:} When $\eta_2\to\infty$, the 
solvent
becomes immobile, and (\ref{mbp0}) reduces to the problem of the motion of a
membrane under the influence of osmosis and surface tension:
\be\label{osmo}
\left.\begin{array}{rclll}
\partial_t c-\frac{\gamma}{\eta_1}\Delta c&=&0&\mbox{ in $\Omega(t)$,} & t > 
0,\\
\alpha H+\gamma c-\eta_3 V_n&=&0&\mbox{ on $\partial \Omega(t)$,} & t > 0,\\
\frac{\gamma}{\eta_1}\partial_nc+cV_n&=&0&\mbox{ on $\partial \Omega(t)$,} & t > 
0.
\end{array}\right\}
\ee
This model has been discussed in one spatial dimension (under the name
``closed osmometer problem'') in  \cite{rub,frh}, in a
radially symmetric setting in \cite{zaal1,zaal2,zaal3} and in higher dimensions
in \cite{lipr1,lipr2}, with the latter reference discussing the two-phase
setting in terms of stability of equilibria.
\item {\bf Fast diffusion:} When $\eta_1 \to 0$, the concentration is forced
on the spatially constant value $c=c(t)=M/|\Omega(t)|$, where $M$ is the total
mass of solute. To our knowledge, the resulting Stokes problem (\ref{mbp0})$_1$,
(\ref{mbp0})$_2$,(\ref{mbp0})$_4$, (\ref{mbp0})$_5$ has not yet attracted any
attention. Starting from (\ref{osmo}), however, one obtains the
surface motion law (\ref{osmo})$_2$, which is just mean curvature flow with a
nonlocal ``braking term'' \cite{meurs}. 
\item {\bf Impermeable membrane:} When $\eta_3\to\infty$, the 
condition $V_n=u\cdot n$ is enforced, i.e. the membrane simply moves according
to the normal component of the velocity field. 
This is the standard kinematic boundary condition for moving liquid surfaces,
also in cases without a membrane. So the resulting problem
\be\label{stokes}
\left.\begin{array}{rclll}
-\eta_2\Delta u+\nabla(q+\gamma c)&=&0&\mbox{ in $\Omega(t)$,} & t > 0,\\
\vdiv u&=&0&\mbox{ in $\Omega(t)$,} & t > 0,\\
\partial_t c-\frac{\gamma}{\eta_1}\Delta c+\nabla c\cdot u&=&0&\mbox{ in 
$\Omega(t)$,} & t > 0,\\
\tau(u,q+\gamma c)n-\alpha Hn&=&0&\mbox{ on $\partial \Omega(t)$,} & t > 0,\\
u\cdot n-V_n&=&0&\mbox{ on $\partial \Omega(t)$,} & t > 0,\\
\partial_nc&=&0&\mbox{ on $\partial \Omega(t)$,} & t > 0
\end{array}\right\}
\ee
describes the free motion of a drop of viscous liquid under the influence of
its own surface tension, combined with a convection-diffusion problem of
solute inside the drop. It is interesting to observe that in the impermeable
case the presence of a solute does not influence the evolution of the domain:
for any smooth evolution $t\mapsto(\Omega(t),c(t))$, the evolution
$t\mapsto(\Omega(t),0)$ satisfies (\ref{stokes}) as well, i.e.
$t\mapsto\Omega(t)$ depends on $\Omega(0)$ only. More precisely, the domain
evolution is given by the Stokes flow problem with surface tension described 
above.  The reason for
this can be understood in the following terms: While the presence of solute 
appears in the Stokes equations (and dynamic boundary condition) only via
a modified pressure term, it is only the velocity field which determines the
domain evolution. 
\end{itemize} 
To normalize all but one of the occurring constants to $1$ we nondimensionalize 
the equations
using the characteristic quantities
\[
L := \frac{\eta_2}{\eta_3}, \quad T := \frac{L^2 \eta_1}{\gamma}, \quad F := 
\frac{L^{N-1} \eta_2}{T}, \quad M := \frac{L F}{\gamma},
\]
for length, time, force, and molarity, respectively. Keeping the same symbols to
 denote dimen-
sionless variables, we rewrite (\ref{mbp0}) in the form
\be\label{nd}
\left.\begin{array}{rclll}
-\Delta u+\nabla(q+c)&=&0&\mbox{ in $\Omega(t)$,} & t > 0,\\
\vdiv u&=&0&\mbox{ in $\Omega(t)$,} & t > 0,\\
\tau(u,q+c)n&=&\kappa Hn&\mbox{ on $\partial \Omega(t)$,} & t > 0,\\
\\
\partial_t c-\Delta c&=& - \nabla c\cdot u&\mbox{ in $\Omega(t)$,} & t > 0,\\
\partial_nc+c(\kappa H+c)&=&0&\mbox{ on $\partial \Omega(t)$,} & t > 0,\\
V_n-\kappa H & = & c+u\cdot n & \mbox{ on $\partial \Omega(t)$,} & t > 0, \\
\\
\Omega(0) & = & \Omega_0, \\
c(0) & = & c_0 & \mbox{ in $\Omega_0$,}  
\end{array}\right\}
\ee
where $\tau$ denotes now the mapping $(u,r) \mapsto \eps(u) - rI$ and $\kappa := \frac{\alpha \eta_1 
\eta_2^{N-1}}{\gamma \eta_3^N}$.

\section{Short-time well-posedness} \label{shwp}

\subsection{Preliminaries and formulation of the main result} \label{pre}

If $U \subset \mathbb{R}^l$ ($l \in \mathbb{N}$) is an  open set, let $BUC(U)$ 
be the Banach space of all bounded and uniformly continuous real valued 
functions on $U$. The space $BUC^k(U)$ contains those elemets of $BUC(U)$ that 
possess bounded and uniformly continuous derivatives up to order $k \in 
\mathbb{N}$. For $k \in \mathbb{N} \cup \{0\}$ and $s \in (0,1)$, $h^{k+s}(U)$ 
denotes the little H\"older space, see \cite{Lu95} for a precise definition and 
basic properties. If $U$ is a domain with sufficiently regular boundary, then 
$h^{k+s}(U)$ is known to be the closure of the smooth functions in the usual 
H\"older space. All spaces are equipped with their natural topologies. As usual, 
function spaces over a manifold are defined by means of a sufficiently smooth 
atlas.\\

First we consider the Stokes equations on a general  domain in the regularity 
scale of little  H\"older spaces:

Let $\beta \in (0,1)$ and $\Xi$ be a bounded $C^{3+\beta}$ domain in 
$\mathbb{R}^N$ (that is a bounded open set whose boundary  possesses a $C^{3 + 
\beta}$-atlas) with exterior unit normal $n=n_{\partial\Xi}$. Let $V_0$ be the 
vector space of
rigid body velocities on $\RRM^N$. Fix a basis $\{\phi_k\}$ in $V_0$ and define the linear operator $\ell\in\cll(L^2(\Xi,\RRM^N),\,V_0)$ by
\[\ell(w)=\ell_\Xi(w)=\sum_k(w,\phi_k)\phi_k,\quad (w,\phi_k):=\int_\Xi
w\cdot\phi_k\,dx.\]
We introduce the spaces
\be\label{spacesXY}
\clx(\Xi):=h^{2+\beta}(\Xi,\RRM^N)\times 
h^{1+\beta}(\Xi),\quad
\cly(\Xi):=h^{\beta}(\Xi,\RRM^N)\times h^{1+\beta}(\Xi)
\times  h^{1+\beta}(\partial \Xi,\mathbb{R}^N).
\ee
\begin{lemma} \label{stokes1}
\begin{itemize}
\item[\rm(i)] The linear operator
\[\Lambda=\Lambda_\Xi\in\cll(\clx(\Xi),\cly(\Xi))\]
given by 
\[\Lambda(u,p):=(-\Delta u+\nabla p+\ell(u),\vdiv u,\,\tau(u,p)n) \]
is an isomorphism.
\item[\rm(ii)] If $(u,p)=\Lambda^{-1}(f,g,h)$ and $(f,g,h)$ satisfies the
solvability conditions
\[\int_\Xi(f-\nabla g)\cdot\phi\,dx-\int_{\partial\Xi}h\cdot\phi\,dS=0\quad
\mbox{ for all $\phi\in V_0$}\]
then $\ell(u)=0$.
\end{itemize} 
\end{lemma}
{\bf Proof:} This result is essentially stated (in large H\"older spaces and
without proof)
in \cite{sol}, Prop. 2, see also p.~645. A proof of a corresponding result in
Sobolev spaces can be found in~\cite[Lemma 2]{gupr}. Basically, the 
proof is
based on the fact that the operator  given by
\[(u,p)\mapsto(-\Delta u+\nabla p,\vdiv u,\,\tau(u,p)n) \]
describes a Douglis-Nirenberg elliptic boundary system that satisfies the
Lopatinskii-Shapiro condition and is therefore a Fredholm operator from
$\clx(\Xi)$ to $\cly(\Xi)$. Its kernel is easily seen to be $V_0\times\{0\}$,
and a discussion of the weak formulation shows that its index is zero and
yields the solvability conditions. Introducing the
auxiliary term $\ell(u)$ accomplishes the reduction to the case of an
isomorphism, cf. \cite{vtr}, Lemma 21.1. \hfill $\blacksquare$
\\ \\
Fix $0 < \beta < \alpha < 1$. We assume that
\begin{itemize}
\item[(I1)] $\Omega_0 \subset \mathbb{R}^N$ is a domain  and $\Gamma_0 :=
\partial \Omega_0$ is a closed compact hypersurface of regularity class
$h^{4+\beta}$;
\item[(I2)] $c_0 \in h^{2+\alpha}(\Omega_0)$  satisfies $\partial_n c_0 + \kappa
H_{\Gamma_0} c_0 + c_0^2 = 0$ on $\Gamma_0$.
\end{itemize}
We use the direct mapping method to transform system (\ref{nd}) into a set of
equations given over a fixed and smooth reference domain. The unknown family of
surfaces $\{ \Gamma(t) \} := \{ \partial \Omega(t) \}$ will be described by a
signed distance function with respect to that surface. In order to carry out
these transformations, we need some preparation:

Given any closed compact hypersurface $\Sigma$  of class $C^2$, let $T_{\delta}
= T_{\delta}(\Sigma)$ be an open tubular neighborhood of $\Sigma$, i.e. the
diffeomorphic image of the mapping 
\[
X_{\Sigma}: \Sigma \times (-\delta,\delta)  \rightarrow \mathbb{R}^N, \qquad
(x,a) \mapsto x + a \, n_{\Sigma} (x),
\]
where $n_{\Sigma}(x)$ is the outer unit normal vector at $x \in \Sigma$ and 
$\delta > 0$ is sufficiently small. It is convenient to decompose the inverse of
$X_{\Sigma}$ into $X_{\Sigma}^{-1}=(P_{\Sigma},\Lambda_{\Sigma})$, where
$P_{\Sigma}(x)$ is the metric projection of a point $x \in T_\delta$ onto
$\Sigma$ and $\Lambda_{\Sigma}$ is the signed distance function with respect to
$\Sigma$. Let
\[
\textrm{Ad}_{\Sigma,\eps}:=\{ \sigma \in C^1(\Sigma); \; \Vert \sigma 
\Vert_{C(\Sigma)} < \eps/5 \} \qquad (\eps > 0)
\]
be the set of admissible boundary perturbations.
If $\eps > 0$ is small enough, then the mapping   $\theta_{\sigma}(x):=x +
\sigma(x)\, n_{\Sigma}(x)$ is for each $\sigma \in
\textrm{Ad}_{\Sigma,\eps}$ a diffeomorphism mapping $\Sigma$ onto
$\Sigma_\sigma:=\theta_\sigma[\Sigma]$.
\\
\\
Due to Theorem 4.2 in \cite{bel} we can fix a number $\delta > 0$ and  a triple
$(\Omega,S_\delta(\Gamma),\rho_0)$ in the following way:
\begin{itemize}
\item $\Omega \subset \Omega_0$ is a domain and $\Gamma:=\partial \Omega$  is a
closed compact real analytic  hypersurface;
\item $S := S_\delta(\Gamma)$ is an open tubular neighborhood of $\Gamma$, $\Gamma_0 \subset S$;
\item $\rho_0 \in h^{4,\beta}(\Gamma) \cap \textrm{Ad}_{\Gamma,\delta}$ and the mapping $\theta_{\rho_0}: \Gamma \rightarrow \Gamma_0$ is a $h^{4,\beta}$ - diffeomorphism. In particular, $\Gamma_0 = \Gamma_{\rho_0}$.
\end{itemize}

From now on let $\delta > 0$, $(\Omega,S,\rho_0)$ be chosen as described above and let $\textrm{Ad}:=\textrm{Ad}_{\Gamma,\delta}$.
\\
\\
Observe that $\sigma[\Gamma] \subset S$ for all $\sigma \in \textrm{Ad}$. Suppose that $\sigma \in \textrm{Ad} \cap h^{m+\gamma}(\Gamma)$ for some $(m,\gamma) \in \mathbb{N} \times (0,1)$. It is not difficult to see that then $\theta_{\sigma} \in h^{m+\gamma}(\Gamma,\RRM^N)$ and $\theta^{-1}_{\sigma} \in h^{m+\gamma}(\Gamma_\sigma,\RRM^N)$. Moreover, given $\sigma \in \textrm{Ad} \cap h^{m+\gamma}(\Gamma)$, the mapping $\theta_\sigma$ extends to a diffeomorphism (the so called Hanzawa - diffeomorphism)
\[
\theta_{\sigma} \in \textrm{Diff}^{m+\gamma}(\mathbb{R}^N,\mathbb{R}^N), \qquad \theta_\sigma |_{\Omega} \in \textrm{Diff}^{m+\gamma}(\Omega,\Omega_{\sigma}) \qquad (\Omega_{\sigma} := \theta_\sigma[\Omega]),
\]
such that we have $\partial \Omega_{\sigma} = \Gamma_\sigma$, see \cite{E04} for details. Note that for $\sigma \in \textrm{Ad}$ the surface $\Gamma_\sigma$ is the zero level set of the function $\varphi_\sigma$ defined by
\[
\varphi_\sigma(x) = \Lambda_{[\Gamma]}(x) - \sigma(P_{[\Gamma]}(x)),
\]
$x \in S$, i.e. $\Gamma_\sigma = \varphi_\sigma^{-1}[\{ 0 \}]$. For later use we set
\[
L_\sigma(x) := |\nabla \varphi_\sigma|(\theta_\sigma(x)).
\]
It can be shown that $L_\sigma > 0$ on $\Gamma$ for all $\sigma \in \textrm{Ad}$. Finally, if $\rho:[0,T] \rightarrow \textrm{Ad}$ is time dependent, we use the notation
\[
\Omega_{\rho,T} := \bigcup_{t \in (0,T)} \{t\} \times \Omega_{\rho(t)} \subset \mathbb{R}^{N+1}.
\]
We are now ready to introduce the notion of a classical solution of (\ref{nd}):
\begin{definition} \label{cs}
Let $c_0$, $\Gamma_0$ satisfy (I1) and (I2), and let
$\mathcal{O}:=h^{4,\beta}(\Gamma) \cap \textrm{Ad}$ inherit  the topology of
$h^{4,\beta}(\Gamma)$. A time-dependent family of domains $\{ \Omega(t); \; t
\in [0,T] \}$, functions $c(t), q(t): \bar\Omega(t) \rightarrow \mathbb{R}$ and
a vector field $u(t): \bar\Omega(t) \rightarrow \mathbb{R}^N$ form a classical
solution of (\ref{nd}) on $[0,T]$, if there exists a function $\rho \in C([0,T],
\mathcal{O}) \cap C^1([0,T],h^{2,\beta}(\Gamma))$ such that letting $\Gamma(t)
:= \partial \Omega(t)$
\begin{itemize}
\item[i)]     $\Omega(t) = \Omega_{\rho(t)}$, $t \in [0,T]$ (thus also
$\Gamma(t) = \Gamma_{\rho(t)}$);
\item[ii)]    $c(\cdot) \circ \theta_{\rho(\cdot)} \in
C([0,T],h^{2+\alpha}(\Omega)) \cap C^1([0,T],h^{\alpha}(\Omega))$;
\item[iii)]   $(u,q)(t) \in h^{3+\beta}(\Omega(t),\mathbb{R}^N) \times
h^{2+\beta}(\Omega(t))$ for $t \in [0,T]$;
\item[iv)]   $t \mapsto (c(t),u(t),q(t),\Omega(t))$  satisfies the equations of
(\ref{nd}) pointwise on $[0,T]$, and, additionally, $l_{\Omega_{\rho(t)}}(u(t))
= 0$ for $t \in [0,T]$.
\end{itemize}
\end{definition} 
Note that ii) in particular implies that $c \in
C^{1,2}(\Omega_{\rho,T},\mathbb{R}) \cap BUC(\Omega_{\rho,T},\mathbb{R})$  and
$c(t) \in h^{2+\alpha}(\Omega_{\rho(t)})$ for $t \in [0,T]$. The main theorem of
this section reads as follows:
\begin{theorem} \label{mt}
Let $c_0$, $\Gamma_0$ satisfy (I1) and (I2). Then  there exists a positive time
$T$ and a unique classical solution $t \mapsto (c(t),u(t),q(t),\Omega(t))$ of
(\ref{nd}) on $[0,T]$.
\end{theorem}

\subsection{Transformation to a fixed interface} \label{trafo}

Given $\sigma \in \textrm{Ad}$,  let $\theta^{*}_{\sigma}$,
$\theta_{*}^{\sigma}$ denote the pull-back and push-forward operators induced by
$\theta_{\sigma}$, i.e. $\theta^{*}_{\sigma} \, f = f \circ \theta_{\sigma}$,
$\theta_{*}^{\sigma} \, g = g \circ \theta_{\sigma}^{-1}$. For functions
$b=b(t,x)$, $\rho=\rho(t,x)$ that depend on time we define $[\theta^{*}_{\rho}
\, b] (t,x) := [\theta^{*}_{\rho(t)} \, b (t,\cdot)](x)$, analogue for
$\theta_{*}^{\rho}$. 

First we consider the transformed Stokes equations.  For $\rho \in \Ad \cap 
h^{3+\beta}(\Gamma)$  observe that (cf.~\cite[Lemma 1]{gupr})
\[\int_{\Gamma_\rho}H_{\Gamma_\rho} n_{\Gamma_{\rho}} \cdot\phi\,dS=0\quad
\mbox{ for all $\phi\in V_0$.}\]
Thus, letting 
\[
\begin{array}{rcl}
H(\rho)  & := & \theta^{*}_{\rho} \, H_{\Gamma_{\rho}}, \\
n(\rho)  & := & \theta^{*}_{\rho} \, n_{\Gamma_{\rho}}, 
\end{array}
\]
in view of Lemma \ref{stokes1} it makes sense to define 
\[s(\rho):=\theta_\rho^\ast v\]
where $(v,q)\in h^{2+\beta}(\Omega_\rho,\RRM^N) \times h^{1+\beta}(\Omega_\rho)$ is the unique solution of
\be\label{stokeseq}
\left.\begin{array}{rcll}
-\Delta v+\nabla q&=&0&\mbox{ in $\Omega_\rho$,}\\
\vdiv v&=&0&\mbox{ in $\Omega_\rho$,}\\
\tau(v,q)n_{\Gamma_\rho}&=& \kappa H_{\Gamma_\rho}n_{\Gamma_\rho}&\mbox{ on
$\Gamma_\rho$,}\\
\ell_{\Omega_\rho}(v)&=&0.
\end{array}\right\}
\ee
This mapping is smooth:
\begin{lemma}\label{stokes2}
We have $[\rho\mapsto s(\rho)]\in C^\infty\big(\Ad\cap
h^{3+\beta}(\Gamma),h^{2+\beta}(\Omega,\RRM^N)\big)$.
\end{lemma}
{\bf Proof:} Recall our notation (\ref{spacesXY}).
For $\rho\in\Ad\cap h^{3+\beta}(\Gamma)$ we have that the
pull-back $\theta_\rho^\ast$
induces isomorphisms from $\clx(\Omega_\rho)$ to $\clx(\Omega)$ and 
from $\cly(\Omega_\rho)$ to $\cly(\Omega)$ which we will denote by the same
symbols. The corresponding inverse will be denoted by $\theta^\rho_\ast$.
Define
\[\Lambda(\rho):=\theta_\rho^\ast\Lambda_{\Omega_\rho}\theta^\rho_\ast,\]
observe that by Lemma \ref{stokes1} 
$\Lambda(\rho)\in \cll_{is}\big(\clx(\Omega),\cly(\Omega)\big)$ and 
\[(s(\rho),\theta_\rho^\ast q)=\Lambda(\rho)^{-1}(0,0,\kappa H(\rho)n(\rho))\]
with $q$ from (\ref{stokeseq}). As 
\[[\rho\mapsto \kappa H(\rho)n(\rho)]\in C^\infty\big(\Ad\cap
h^{3+\beta}(\Gamma),h^{1+\beta}(\Omega,\RRM^N)\big),\]
it remains to show that
\be\label{lamsmooth} [\rho\mapsto \Lambda(\rho)]\in
C^\infty\big(\Ad\cap
h^{3+\beta}(\Gamma),\cll(\clx(\Omega),\cly(\Omega)\big).
\ee
The result follows then from the fact that taking the inverse of an isomorphism
is a smooth operation. To show (\ref{lamsmooth}) it is sufficient to explicitly
carry out the transformations of the differential and integral operators
involved. For example, for a first-order partial derivative
$\partial_i$ we have
\[\theta_\rho^\ast\partial_i\theta^\rho_\ast=a^j_i\partial_j,\]
where $a^j_i=a^j_i(\rho)\in h^{2+\beta}(\bar \Omega)$ is an element of the matrix
$(D\theta_\rho)^{-1}$ and therefore depends smoothly on $\rho\in\Ad\cap
h^{3+\beta}(\Gamma)$. Similarly,
\[\theta_\rho^\ast(\ell_{\Omega_\rho}(\theta^\rho_\ast v)) =
\sum_i\int_\Omega v\cdot(\psi_i\circ\theta_\rho)\det
D\theta_\rho\,dx\,(\psi_i\circ\theta_\rho),\]
which is smooth in $\rho$ as well. \hfill $\blacksquare$
\\ \\
In order to transform (\ref{nd}) to $\Omega$ for suitable $\rho$ we introduce the operators $\mathcal{A}(\rho)$, $\mathcal{B}(\rho)$, $\mathcal{K}(\rho)$ by 
\[
\begin{array}{rcl}
\mathcal{A}(\rho)\xi & := & \theta^{*}_{\rho} (\Delta (\theta_{*}^{\rho} \xi)); \\
\mathcal{B}(\rho)\xi & := & \theta^{*}_{\rho} (\nabla (\theta_{*}^{\rho} \xi)) \cdot n(\rho); \\
\mathcal{K}(\rho)\xi & := & \theta^{*}_{\rho} (\nabla (\theta_{*}^{\rho} \xi)).
\end{array}
\]
The transformed problem reads then
\be\label{ndtr}
\left.\begin{array}{rcll}
\partial_t \xi - \cla(\rho) \xi & = & R(\xi,\rho) - \mathcal{K}(\rho)\xi \cdot s(\rho) & \mbox{ in $\Omega$,}\\
\clb(\rho) \xi & = & -\kappa \xi H(\rho) - \xi^2 & \mbox{ on $\Gamma$,}\\
\partial_t \rho - \kappa L_\rho H(\rho) & = & L_\rho(\xi + s(\rho) \cdot n(\rho))
&\mbox{ on $\Gamma$,}\\
\xi(0)&=&\xi_0,\\
\rho(0)&=&\rho_0,
\end{array}\right\}
\ee
where $\xi_0:=\theta^{*}_{\rho_0} c_0$. The term $R$ arises from the transformation of the time derivative and is determined by
\[
R(z,\sigma)(y) = r_0(L_\sigma [\kappa H(\sigma) + z + s(\sigma) \cdot n(\sigma)],B_{\mu}(\sigma)z)(y), \qquad y \in \Omega,
\]
where $z \in C^1(\bar{\Omega})$,  $\sigma \in \textrm{Ad} \cap C^2(\Gamma)$ and
\begin{equation}
\label{BLM}
r_0(h,k)(y):=
\left\{ \begin{array}{lcl}
\chi(\Lambda_\Gamma(y)) \cdot h(P_{\Gamma}(y)) \cdot k(y), & \textrm{if} & y \in \Omega \cap S \\
0,                                            & \textrm{if} & y \in \Omega \setminus (\Omega \cap S), \end{array} \right.
\end{equation}
\[
B_{\mu}(\sigma)z (y) =  \theta^{*}_{\sigma} \, \nabla (\theta_{*}^{\sigma} z) (y) \cdot (n_{\Gamma} \circ P_{\Gamma})(y) , \qquad y \in S
\]
($\chi$ being a suitable cut-off function, cf. \cite{E04}, \cite{K07}). The derivation of (\ref{ndtr}) is a straightforward calculation \cite{E04}, \cite{K07}.
\\
\\
Note: If $(\xi,\rho)$ is a sufficiently regular solution of (\ref{ndtr}), then $(\theta^\rho_* \xi, \theta^\rho_* s(\rho), q, \Gamma_\rho)$ is a classical solution of (\ref{nd}), where $q := \theta^\rho_* (P_2 \Lambda(\rho)^{-1}(0,0,\kappa H(\rho)n(\rho)) - \xi)$ ($P_2$ denoting the projection on the second component).
\\
\\
It turns out that the local well-posedness of system (\ref{ndtr}) can be proved almost in the same fashion as outlined in sections 2 and 3 in \cite{lipr1}, so we recall the abstract setting from there: If  $T > 0$ is given and $J_T:=[0,T]$, let
\[\begin{array}{rcl}
E_0 & := & h^{\alpha}(\Omega)\times h^{2+\beta}(\Gamma), \\
E_1 & := & h^{2+\alpha}(\Omega)\times h^{4+\beta}(\Gamma),\\
\bbe_0(J_T) & := & BUC(J_T,E_0), \\
\bbe_1(J_T) & := & BUC^1(J_T,E_0)\cap BUC(J_T,E_1),\\
\bbf(J_T)   & := & BUC(J_T,h^{1+\alpha}(\Gamma))\cap h^{(1+\alpha)/2}(J_T,C(\Gamma)).
\end{array}\]
To economize notation we drop the $T$ - dependence, i.e. write $\bbe_1$ instead of $\bbe_1(J_T)$ etc. and define the sets
\[\widetilde{\mbox{Ad}}=\{(\nu,\psi)\in E_1\,|\,\psi\in\mbox{Ad}\},
\quad
\widehat{\mbox{Ad}}=\{w\in\bbe_1\,|\,w(t)\in\widetilde{\mbox{Ad}},\,t\in[0,T]\},\]
which are open subsets of $E_1$ and $\bbe_1$, respectively. Our goal is to write system (\ref{ndtr}) as a single operator equation. For this we recall the splitting
\[
H(\rho) = P(\rho)\rho + Q(\rho)
\]
as described for example in \cite{ES97a}. More precisely, $P(\rho)$ can be chosen to be a second order uniformly  elliptic operator acting as an isomorphism in various scales of function spaces and depending smoothly on $\rho$. The mapping $Q$ contains only lower order terms. Precise mapping properties of $P(\cdot)$ and $Q$ are given for example in \cite{E04}, \cite{ES97a}, \cite{K07}. Let  
\bea
\bba(w)(t)&=&\left(\begin{array}{cc}
\cla(\rho(t))&0\\
0&\kappa L_{\rho(t)} P(\rho(t))
\end{array}\right),\\
\tilde\bbb(\nu,\psi)(\zeta,\chi)&=&\clb(\psi)\zeta,\\
\big(\bbb(w)(v,\sigma)\big)(t)&=&\tilde\bbb(w(t))(v(t),\sigma(t)),\\
\bbl(w)&=&(\partial_t-\bba(w),\bbb(w),\gamma_t),
\eea
where $w=(\xi,\rho)\in\widehat{\mbox{Ad}}$, $(v,\sigma)\in\bbe_1$, $(\nu,\psi) \in \widetilde{\mbox{Ad}}$, $(\zeta,\chi)\in E_1$ and $\gamma_t\in\cll(\bbe_1,E_1)$  denotes the time trace map $w\mapsto w(0)$. We have 
\bea
\bba&\in& C^\infty\big(\widehat{\mbox{Ad}},\cll(\bbe_1,\bbe_0)\big),\\
\tilde\bbb&\in& C^\infty\big(\widetilde{\mbox{Ad}},\cll(E_1,h^{1+\alpha}(\Gamma))\big),\\
\bbb&\in& C^\infty\big(\widehat{\mbox{Ad}},\cll(\bbe_1,\bbf)\big),\\
\bbl&\in& C^\infty\big(\widehat{\mbox{Ad}},
\cll(\bbe_1,\bbe_0\times\bbf\times E_1)\big),
\eea
cf. \cite{lipr1}. Let $w_0=(\xi_0,\rho_0)$. For given, fixed $M > \Vert w_0 \Vert_{E_1}$ we define the closed set
\[ \clc = \clc(M,T) := \{w\in\bbe_1\,|,w(0)=w_0,\|w\|_{\bbe_1}\leq M\}\]
and introduce the subspace $\bbz \subset \bbe_0\times\bbf\times E_1$ by 
\[\bbz=\{(f,g,h)\in\bbe_0\times\bbf\times E_1\,|\,\gamma_tg=\tilde\bbb(w_0) h\}.\]
The following lemma collects some facts shown in \cite{lipr1} (Lemmas 3.1 - 3.5). The symbol  $\cll_{is}$ stands for the set of topological isomorphisms.
\begin{lemma} \label{prop} 
Let $M > \Vert w_0 \Vert_{E_1}$. There is a $T^\ast=T^\ast(M,w_0)$ and a $C=C(w_0)$ such that if $T\in(0,T^\ast]$ then $\clc\subset\widehat{\mbox{Ad}}$, $\bbl(\clc)\subset\cll_{is}(\bbe_1,\bbz)$ and
\[\|\bbl(w)^{-1}\|_{\cll(\bbz,\bbe_1)}\leq C,\qquad w\in\clc.\]
\end{lemma}
Thus, our problem can be reformulated as
\be\label{prob}
\bbl(w)w=F(w):=(\clr(w),\clg(w),w_0), \quad w\in\clc,
\ee
where
\bea
\clr(w)(t)&=&
\left(\begin{array}{c}R(w(t)) - \mathcal{K}(\rho(t))\xi(t) \cdot s(\rho(t))\\L_{\rho(t)}(\kappa Q(\rho(t)) + \xi(t) + s(\rho(t)) \cdot n(\rho(t)))
\end{array}\right),\\
\clg(w)(t)&=&-\kappa \xi(t)H(\rho(t))-\xi(t)^2,
\eea
$w=(\xi,\rho)$. In view of Lemma \ref{stokes2} it is easily checked that 
\[F\in C^\infty (\widehat{\mbox{Ad}},\bbe_0\times\bbf\times E_1),\]
cf. \cite{lipr1}. 
\begin{lemma} (Quasilinear character)\label{qualin}
Let $\eps>0$ and $M > \Vert w_0 \Vert_{E_1}$ be given. There is a $T^\ast=T^\ast(\eps,M,w_0)$ such that
if $T\in(0,T^\ast]$, $w_1,w_2\in\clc$, then
\begin{eqnarray}\label{quaestL}
\|\bbl(w_1)-\bbl(w_2)\|_{\cll(\bbe_1,\bbz)}&\leq&\eps\|w_1-w_2\|_{\bbe_1};\\
\|F(w_1)-F(w_2)\|_{\bbz}&\leq&\eps\|w_1-w_2\|_{\bbe_1}.\label{quaestF}
\end{eqnarray}
\end{lemma}
{\bf Proof:} The estimate (\ref{quaestL}) has been proven in \cite{lipr1}. Using Lemma \ref{stokes2} and the facts that 
\[
\begin{array}{l}
\mathcal{K} \in C^{\infty}\big(\mbox{Ad}  \cap 
h^{3+\beta}(\Gamma),\mathcal{L}(h^{1+\alpha}(\Omega),h^{\alpha}(\Omega,\mathbb{R
}^N))\big), \\
R \in C^{\infty}\big(\mbox{Ad} \cap h^{3+\beta}(\Gamma) \times h^{1+\alpha}(\Omega),h^{\alpha}(\Omega)\big), \\
n \in C^{\infty}\big(\mbox{Ad} \cap h^{3+\beta}(\Gamma),h^{2+\beta}(\Gamma,\mathbb{R}^N)\big),
\end{array}
\]
the estimate (\ref{quaestF}) results analogously to A.17 and A.18 in 
\cite{lipr1}.  Observe in this connection that pointwise scalar multiplication 
canonically induces a bounded and bilinear mapping
\[
h^{m+\gamma}(M,\mathbb{R}^l) \times h^{\tilde m+\tilde \gamma}(M,\mathbb{R}^l)  
\rightarrow h^{m+\gamma}(M,\mathbb{R}), \qquad  M \in \{ \Gamma, \Omega \}, 
\]
where $m,\tilde m,l \in \mathbb{N} \cup \{0\}$, $l \geq 1$, $\gamma,  \tilde 
\gamma \in (0,1)$ and $\tilde m + \tilde \gamma \geq m + \gamma$. \hfill 
$\blacksquare$ 
\begin{theorem} (Short-time wellposedness) \label{exist}
Let $w_0=(\xi_0,\rho_0)\in h^{2+\alpha}(\bar\Omega)
\times(h^{4+\beta}(\Gamma)\cap\mbox{Ad})$ be as specified above. Then there  are 
constants $M,T^\ast>0$ such that {\rm(\ref{prob})}, or, equivalently, 
{\rm(\ref{ndtr})} has precisely one solution in $\clc$ for any $T\in(0,T^\ast]$.
\end{theorem}
{\bf Proof:} Observe that $(\xi_0,\rho_0)$ satisfy the compatibility condition
\[\clb(\rho_0)\xi_0=-\kappa \xi_0 H(\rho_0)-\xi_0^2\]
because of (I2). Moreover, due to Lemma \ref{prop} we can rewrite (\ref{prob}) 
as  a fixed point equation
\be\label{fpoint}
w=\Phi(w):=\bbl(w)^{-1} F(w),\qquad w\in \mathcal{C}.
\ee
Thus, the assertion is an immediate consequence of Lemmas \ref{prop}, 
\ref{qualin}  and the obvious modifications of Lemmas 3.5, 3.6 in \cite{lipr1}. 
\hfill $\blacksquare$
\\ \\
The statement of Theorem \ref{mt} is a direct consequence of Theorem \ref{exist} 
 and the observation that for a $C^{4+\beta}$ domain $\Xi$ we have  
\[\Lambda_\Xi \in \mathcal{L}_{\mbox{is}} \big(h^{3+\beta}(\Xi,\RRM^N)\times 
h^{2+\beta}(\Xi), h^{1+\beta}(\Xi,\RRM^N)\times h^{2+\beta}(\Xi) \times  
h^{2+\beta}(\partial \Xi,\mathbb{R}^N)\big).\]

\section{Conclusion}

Our modelling approach consisted essentially in ``adding up building blocks'' 
from problems with a well-known variational structure, namely, diffusion and 
Stokes flow with surface tension, for the free energy functional as well as for 
the dissipation.  Even though our evolution  can be interpreted 
as a gradient flow, the character and structural properties of the resulting 
nonlinear problem are (to us) not a priori obvious. In particular,  
even with  the same state space and dissipation functional, different energy 
functionals may lead to both parabolic and hyperbolic evolutions. An example 
for this is given by the space of probability measures on the real axis with 
the Wasserstein metric, where, as is well-known by now, the (generalized) 
gradient flow with respect to the entropy functional  is the heat flow, while a 
certain class of autocorrelation functionals gives rise to a nonlocal 
hyperbolic evolution related to Burger's equation \cite{BCFP}. In this respect, 
the challenging problem arises to find direct connections between structural 
conditions on the energy and dissipation functionals on one hand and the type 
or other properties of the corresponding evolution on the other.  At the 
moment, we feel unable to even give reasonably general nontrivial conjectures 
on this.

It turns out that  in our case the resulting evolution is parabolic in the 
following sense: The associated linear homogeneous evolution is described by an 
analytic semigroup of operators. In turn, the theory of these semigroups provides 
the means to prove optimal regularity results for the corresponding linear, 
nonhomogeneous evolution equations that arise from linearizing the original problem. 
For a more precise discussion of this, we refer to \cite{LaQPP} or \cite{Lu89}.

In our problem, we have to consider a coupled evolution for a pair of 
functions, one of them given inside the reference domain (with boundary 
conditions) and the other on its boundary. The generator of the corresponding 
semigroups is diagonal in highest order, so that known results on the 
``components'' can be applied, including a crucial optimal regularity result 
(Theorem 1.4) from \cite{Lu89}. Technically, this is the basis for Lemma 
\ref{prop} in the present paper. 

Furthermore, it is important for our analysis that the nonlocal 
solution operator of the Stokes equations only enters in a lower order term  
(since the pseudodifferential operator mapping the Neumann normal stress 
boundary data to the  Dirichlet boundary data is of order $-1$, cf. Lemma 
\ref{stokes2}), and therefore  does not occur in 
the leading linear operator $\mathbb{L}(w_0)$ defined in  Section \ref{trafo}.

\newpage
\subsection*{Acknowledgements}
The authors are indebted to M.M. Zaal for interesting discussions. The first
author was supported by DFG SPP 1506 ``Transport processes at fluidic
interfaces''.

\end{document}